\documentclass[hoptionsi]{smfart}

\usepackage[francais,english]{babel}
\usepackage{smfthm}

\usepackage[all]{xy}

\usepackage{amsmath}
\usepackage{amsxtra}
\usepackage{amscd}
\usepackage{amsthm}
\usepackage{amsfonts}
\usepackage{amssymb}
\usepackage{eucal}

\usepackage{graphicx,epsfig}

\newcommand{\Glie}{\mathfrak{g}}

\newcommand{\U}{\mathcal{U}}

\newcommand{\ZZ}{\mathbb{Z}}
\newcommand{\CC}{\mathbb{C}}

\renewcommand{\theequation}{\thesection.\arabic{equation}}
\def\theenumi{\roman{enumi}}
\def\labelenumi{(\theenumi)}

\newcommand{\nc}{\newcommand}
\nc{\on}{\operatorname}
\nc{\la}{\lambda}
\nc{\wh}{\widehat}

\author[David Hernandez]{David Hernandez}

\address{Sorbonne Paris Cit\'e, Univ Paris Diderot, CNRS Institut de
  Math\'ematiques de Jussieu-Paris Rive Gauche UMR 7586,
B\^atiment Sophie Germain, Case 7012,
75205 Paris Cedex 13, France.}

\email{david.hernandez@imj-prg.fr}

\urladdr{http://webusers.imj-prg.fr/~david.hernandez/}

\title[Spectre quantique]{Spectre des syst\`emes int\'egrables quantiques et repr\'esentations lin\'eaires}

\begin{document}


\frontmatter

\begin{abstract} La structure des valeurs propres d'un syst\`eme quantique, c'est-\`a-dire de son spectre, est essentielle \`a sa compr\'ehension. Dans un c\'el\`ebre article dat\'e de 1971, Baxter a calcul\'e ces valeurs propres pour le mod\`ele \og de la glace\fg. Il a montr\'e qu'elles ont une forme remarquable et r\'eguli\`ere faisant intervenir des polyn\^omes. Dans les ann\'ees 1980-1990, il a \'et\'e conjectur\'e que de tels polyn\^omes permettent de d\'ecrire le spectre de nombreux syst\`emes quantiques plus g\'en\'eraux. Nous allons voir comment, en adoptant le point de vue math\'ematique de la th\'eorie des repr\'esentations, ces polyn\^omes (de Baxter) apparaissent naturellement. Ce r\'esultat nous a permis de d\'emontrer en 2013 la conjecture g\'en\'erale.
\end{abstract}

\maketitle

\mainmatter

\section{Syst\`emes int\'egrables quantiques}

Le {\bf mod\`ele \`a 6 sommets} est un c\'el\`ebre mod\`ele de {\bf physique statistique} introduit par Pauling en 1935, qui permet notamment 
de d\'ecrire le cristal de la glace  (voir \cite{BaxterBook}). Il est r\'ealis\'e sur un r\'eseau dont chaque
sommet est reli\'e \`a 4 autres sommets. Un \'etat du syst\`eme est une orientation
des ar\^etes telle qu'\`a chaque sommet arrivent exactement 2 fl\^eches (Figure 1). 
Les fl\^eches repr\'esentent l'orientation des mol\'ecules d'eau du cristal
les unes par rapport aux autres. 
Il y a 6 configurations
possibles \`a chaque sommet (Figure 2), ce qui justifie l'appellation de ce mod\`ele.
~\begin{center}
\begin{figure}\label{orientation}
 \hspace{4.5cm}  \epsfig{file=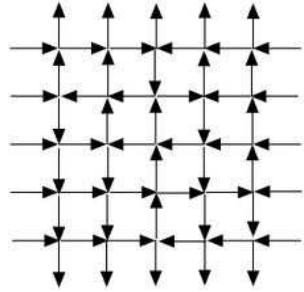,width=0.3
   \linewidth}
      \caption{Une orientation d'un r\'eseau (mod\`ele \`a 6 sommets).}
\end{figure}
\end{center}
\begin{center}
\begin{figure}\label{6vertex}
  \hspace{1cm} \epsfig{file=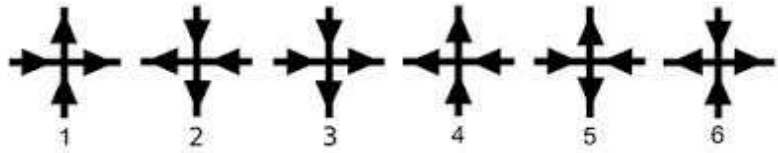,width=0.8\linewidth}
      \caption{6 configurations possibles \`a chaque sommet.}
\end{figure}
\end{center}

L'\'etude du mod\`ele de la glace est fortement li\'ee \`a celle d'un autre mod\`ele, cette fois-ci en {\bf physique statistique quantique}, appel\'e {\bf mod\`ele $XXZ$} de Spin $1/2$, dit de Heisenberg quantique (1928). Il s'agit d'une variante en physique quantique du mod\`ele d'Ising (1925) (voir \cite{JM}), qui mod\'elise des cha\^ines de spins magn\'etiques quantiques
ayant deux \'etats classiques, haut ou bas (Figure 3).
\begin{center}
\begin{figure}\label{spin}
 \hspace{4.5cm}  \epsfig{file=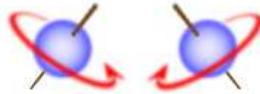,width=0.3
   \linewidth
   }
      \caption{Etats d'un Spin $1/2$ (haut ou bas).}
\end{figure}
\end{center}

Ces deux mod\`eles, mod\`ele \`a 6 sommets et mod\`ele $XXZ$, figurent parmi les plus \'etudi\'es
en physiques statistique et quantique. Les structures math\'ematiques qui les sous-tendent sont tr\`es proches.
En d\'epit de leur formulation assez \'el\'ementaire,
ils sont extr\^emement riches et leur analyse a une tr\`es longue histoire.

En physique statistique (quantique), le comportement du syst\`eme est contr\^ol\'e par la {\bf fonction de partition} $\mathcal{Z}$ \footnote{En physique statistique, la fonction de partition s'exprime comme la somme $\sum_j \text{exp}(-E_j/(k_BT))$ sur tous les \'etats $j$ du syst\`eme, o\`u $E_j$ est l'\'energie de l'\'etat $j$, $T$ est la temp\'erature du syst\`eme et $k_B$ la constante de Boltzmann. En physique quantique, la somme est remplac\'ee par une trace $\text{Tr}_W(\text{exp}(-E/(k_BT)))$ o\`u $E$ est l'op\'erateur \og hamiltonien\fg{}  qui agit sur l'espace $W$ des \'etats quantiques du syst\`eme.}, qui permet d'obtenir les grandeurs mesurables\footnote{Une grandeur mesurable $Q$ est obtenue comme moyenne pond\'er\'ee sur les \'etats $\frac{\sum_j  \text{exp}(-E_j/(k_BT)) Q_j}{\mathcal{Z}}$ des valeurs $Q_j$ sur chaque \'etat $j$.}. Cette fonction
$\mathcal{Z}$ est tr\`es difficile \`a calculer en g\'en\'eral. La m\'ethode de la {\bf matrice de transfert}
est un proc\'ed\'e pour tenter de la d\'eterminer : il s'agit d'\'ecrire $\mathcal{Z}$ comme trace
d'un op\'erateur $\mathcal{T}$ (la matrice de transfert) agissant sur l'{\bf espace des \'etats} $W$ : 
$$\mathcal{Z} = \on{Tr}_W (\mathcal{T}^M).$$
Ici $M$ est un entier associ\'e \`a la taille du r\'eseau du mod\`ele.
Ainsi, pour trouver $\mathcal{Z}$, il suffit d'obtenir les valeurs propres $\lambda_j$ de $\mathcal{T}$ :
$$\mathcal{Z} = \sum_j \lambda_j^M.$$
Le spectre $\{ \lambda_j \}_j$ de $\mathcal{T}$ est appel\'e le {\bf spectre du syst\`eme quantique}.

Dans un c\'el\`ebre article s\'eminal de 1971, inspir\'e notamment par les travaux de Bethe (1931), Baxter \cite{Baxter} a compl\`etement r\'esolu ce probl\`eme\footnote{Baxter a introduit la m\'ethode puissante des \og Q-op\'erateurs\fg{} qui lui a \'egalement permis de r\'esoudre le mod\`ele \og \`a 8 sommets\fg, plus complexe. Le mod\`ele \`a 6 sommets avait aussi \'et\'e r\'esolu par d'autres m\'ethodes, notamment dans les travaux de Lieb et Sutherland (1967).}. Gr\^ace \`a une \'etude
tr\`es pr\'ecise il a notamment montr\'e
que les valeurs propres $\lambda_j$ de $\mathcal{T}$ ont une structure tout \`a fait remarquable : elles
s'expriment sous la forme
\begin{equation}    \label{relB}
\lambda_j = A(z) \frac{Q_j(zq^2)}{Q_j(z)} + D(z)  \frac{Q_j(zq^{-2})}{Q_j(z)},
\end{equation}
o\`u $z,q\in\CC^*$ sont des param\`etres du mod\`ele (respectivement spectral et quantique),
$A(z)$ et $D(z)$ sont des fonctions \og universelles\fg{} (au sens o\`u elles ne d\'ependent pas de la valeur propre
$\lambda_j$). La fonction $Q_j(z)$ d\'epend de la valeur propre, mais c'est un polyn\^ome.
La relation (\ref{relB}) est la fameuse {\bf relation de Baxter} (ou \og relation $TQ$ de Baxter\fg).
Les polyn\^omes $Q_j$ sont appel\'es {\bf polyn\^omes de Baxter}.

\medskip

En r\'esultent alors naturellement les questions suivantes :

- Y a-t-il une explication pour l'existence de la relation de Baxter  ?

- Une expression analogue avec des polyn\^omes permet-elle de d\'ecrire le spectre d'autres syst\`emes quantiques ?

\medskip

Une conjecture formul\'ee en 1998 par Frenkel-Reshetikhin \cite{Fre} affirme que la deuxi\`eme question doit avoir une r\'eponse positive. Comme on ne peut esp\'erer effectuer en g\'en\'eral le calcul d\'etaill\'e de Baxter qui est connu pour le mod\`ele $XXZ$, c'est en r\'epondant \`a la premi\`ere question que nous pouvons d\'emontrer cette conjecture. Pour ce faire, \'etudions les structures math\'ematiques, alg\'ebriques, sous-jacentes \`a la th\'eorie.

\section{Groupes quantiques et leurs repr\'esentations}

Les {\bf groupes quantiques} sont graduellement apparus au cours des ann\'ees 1970, en particulier dans les travaux de l'\'ecole de Leningrad, comme le cadre naturel math\'ematique pour \'etudier les matrices de transfert. Drinfeld \cite{Dri} et Jimbo \cite{J} ont ind\'ependamment d\'ecouvert une formulation alg\'ebrique uniforme sous forme d'{\bf alg\`ebres de Hopf}. Il s'agit d'un des r\'esultats cit\'es pour la m\'edaille Fields de Drinfeld en 1990.

Pour introduire les groupes quantiques de Drinfeld-Jimbo, consid\'erons d'abord un objet tr\`es classique, une {\bf alg\`ebre de Lie} (simple) complexe de dimension finie. Il s'agit d'un espace vectoriel de dimension finie $\mathfrak{g}$ muni d'un {\bf crochet de Lie}, c'est-\`a-dire d'une application bilin\'eaire antisym\'etrique 
$$[,]:\mathfrak{g}\times \mathfrak{g}\rightarrow \mathfrak{g}$$ 
satisfaisant la formule de Jacobi
$$[x,[y,z]] + [y,[z,x]] + [z,[x,y]] = 0\text{ pour tous }x,y,z\in\mathfrak{g}.$$
L'exemple le plus simple, mais n\'eanmoins non trivial car il correspond au mod\`ele $XXZ$, est celui de l'alg\`ebre de Lie $\mathfrak{g} = sl_2$ : c'est l'espace des matrices complexes $2\times 2$ de trace nulle, muni du crochet naturel 
$$[M,N] = MN - NM$$ 
pour lequel il est clairement stable. Pour les g\'en\'erateurs lin\'eaires
$$E = \begin{pmatrix}0&1\\0&0\end{pmatrix}\text{ , }F = \begin{pmatrix}0&0\\1&0\end{pmatrix}
\text{ , }H = \begin{pmatrix}1&0\\0&-1\end{pmatrix},$$
on a par exemple la relation
\begin{equation}\label{crochet}[E,F] = H.\end{equation}

Ces alg\`ebres de Lie ont des analogues naturelles de dimension infinie, les {\bf alg\`ebres de lacets} 
$$\hat{\Glie} = \Glie \otimes \CC[t^{\pm 1}],$$
avec le crochet de Lie d\'efini par
$$[x\otimes f(t),y\otimes g(t)] = [x,y]\otimes (fg)(t)\text{ pour $x,y\in\Glie$ et $f(t),g(t)\in \CC[t^{\pm 1}]$},$$ 
ce qui revient \`a remplacer le corps $\CC$ par l'anneau des polyn\^omes de Laurent complexes 
$$\CC[t^{\pm 1}] = \left\{\sum_{N\leq i\leq M} a_i t^i| N,M\in\ZZ, a_i\in\CC\right\}.$$ 
Ces alg\`ebres sont des quotients d'{\bf alg\`ebres de Kac-Moody affines}, qui ont des propri\'et\'es alg\'ebriques semblables \`a celles des alg\`ebres de Lie simples de dimension finie (notamment une pr\'esentation analogue \`a celle de Serre pour $\Glie$, comme l'ont montr\'e Kac (1968) et Moody (1969), voir \cite{ka}). Elles ont \'et\'e \'etudi\'ees intensivement pour leurs diverses applications en math\'ematiques et physique math\'ematique.

Maintenant, pour \'etudier les syst\`emes quantiques qui nous int\'eressent, ces alg\`ebres de Lie classiques doivent \^etre {\bf quantifi\'ees}, c'est-\`a-dire d\'eform\'ees en tenant compte du param\`etre quantique 
$$q = \text{exp}(h)\in\CC^*,$$ 
o\`u $h$ est un analogue de la grandeur de Planck ($q$ sera bien identifi\'e au param\`etre quantique de la relation (\ref{relB})). 
On retrouve les structures classiques pour $h\rightarrow 0$, donc $q\rightarrow 1$. {\it On supposera dans la suite que $q$ n'est pas une racine de l'unit\'e.}

Bien qu'une telle quantification des alg\`ebres de Lie $\Glie$ ou $\hat{\Glie}$ elles-m\^emes ne soit pas connue, Drinfeld et Jimbo ont d\'ecouvert qu'il existe une quantification naturelle de leurs {\bf alg\`ebres enveloppantes} respectives $\mathcal{U}(\Glie)$ et $\mathcal{U}(\hat{\Glie})$ (alg\`ebres universelles d\'efinies \`a partir des alg\`ebres de Lie, par exemple en rempla\c cant dans la pr\'esentation de Serre les crochets $[x,y]$ par des expressions alg\'ebriques $xy - yx$ dans l'alg\`ebre). On obtient alors les groupes quantiques $\mathcal{U}_q(\Glie)$, $\mathcal{U}_q(\hat{\Glie})$ qui d\'ependent\footnote{Elles peuvent \^etre d\'efinies comme des alg\`ebres sur $\CC[[h]]$.} du param\`etre quantique $q$, voir \cite{CP}. 

Par exemple dans $\mathcal{U}_q(sl_2)$ la relation (\ref{crochet}) devient
$$[E,F] = \frac{e^{hH} - e^{-hH}}{q - q^{-1}}$$
qui tend bien vers $H$ quand $h$ tend vers $0$.

Le cas des {\bf alg\`ebres affines quantiques} $\mathcal{U}_q(\hat{\Glie})$ est particuli\`erement remarquable car Drinfeld \cite{Dri2} a d\'emontr\'e\footnote{La preuve a \'et\'e pr\'ecis\'ee par la suite par Beck puis par Damiani.} qu'elles peuvent non seulement \^etre obtenues comme quantification de $\mathcal{U}(\hat{\Glie})$, mais \'egalement, par un autre proc\'ed\'e, comme {\bf affinisation} du groupe quantique $\mathcal{U}_q(\Glie)$. C'est la {\bf r\'ealisation de Drinfeld} des alg\`ebres affines quantiques. Ceci peut \^etre \'enonc\'e dans le diagramme \og commutatif\fg{} suivant :
$$\xymatrix{ &\hat{\Glie}\ar@{-->}[dr]^{\text{Quantification}}&   
\\ \Glie \ar@{-->}[ur]^{\text{Affinisation}}\ar@{-->}[dr]_{\text{Quantification}}& &  \U_q(\hat{\Glie})
\\ & \U_q(\Glie)\ar@{-->}[ur]_{\text{Affinisation quantique}}& }$$
Ce th\'eor\`eme, qui revient \`a donner deux pr\'esentations isomorphes de $\mathcal{U}_q(\hat{\Glie})$, est un analogue quantique du th\'eor\`eme classique de Kac et Moody. Il s'agit d'une bonne indication de l'importance de ces alg\`ebres d'un point de vue alg\'ebrique.

Les alg\`ebres affines quantiques $\U_q(\hat{\Glie})$ ont en fait une structure beaucoup plus riche, ce sont des alg\`ebres de Hopf. Elles sont notamment munies d'une {\bf comultiplication} (qui est une op\'eration duale de la multiplication), c'est-\`a-dire d'un morphisme d'alg\`ebre
\begin{equation}\label{coproduit}\Delta : \U_q(\hat{\Glie}) \rightarrow \U_q(\hat{\Glie})\otimes \U_q(\hat{\Glie}).\end{equation}
Mais surtout, $\U_q(\hat{\Glie})$ poss\`ede une {\bf $R$-matrice universelle}, c'est-\`a-dire un \'el\'ement (canonique) dans le carr\'e tensoriel\footnote{En fait, dans une l\'eg\`ere compl\'etion du carr\'e tensoriel.}
$$\mathcal{R}(z) \in (\U_q(\hat{\Glie})\otimes \U_q(\hat{\Glie}))[[z]]$$
qui est notamment une solution de l'{\bf \'equation de Yang-Baxter quantique} :
$$\mathcal{R}_{12}(z)\mathcal{R}_{13}(zw)\mathcal{R}_{23}(w) = \mathcal{R}_{23}(w)
\mathcal{R}_{13}(zw) \mathcal{R}_{12}(z).$$
Les param\`etres formels $z$, $w$ sont appel\'es {\bf param\`etres spectraux}. Cette \'equation est \`a valeurs dans le cube tensoriel 
$$(\U_q(\hat{\Glie}))^{\otimes 3}[[z, w]].$$
Les indices dans les facteurs indiquent l'emplacement des termes de la $R$-matrice universelle : 
$$\mathcal{R}_{12}(z) = \mathcal{R}(z)\otimes 1\text{ , }\mathcal{R}_{23}(z) = 1 \otimes \mathcal{R}(z)...$$
Il s'agit d'une \'equation hautement non triviale, li\'ee aux mouvements de tresses. En effet, dans la figure 4 on retrouve l'\'equation en lisant de bas en haut et en multipliant
par un facteur $\mathcal{R}_{\alpha\beta}$ d'indice $(\alpha,\beta)$ lorsque le brin $\alpha$ croise le brin $\beta$. C'est pour cette raison que la th\'eorie des repr\'esentations des groupes quantiques permet de construire des invariants en topologie de basse dimension (notamment les polyn\^omes de Jones des n\oe uds). Il s'agit historiquement, avec la construction par Lusztig et Kashiwara de bases canoniques de repr\'esentations des alg\`ebres de Lie classiques, d'un des premiers grands succ\`es de la th\'eorie des groupes quantiques. Nous n'aborderons pas ces sujets ici pour nous concentrer sur les applications aux syst\`emes quantiques.
 \begin{center}
\begin{figure}\label{tresse}
 \hspace{4cm}  
 {\epsfig{file=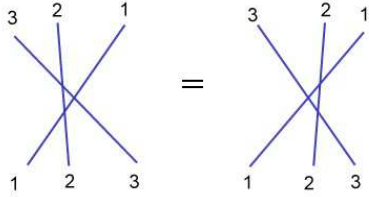,width=.4\linewidth}
}
      \caption{Equation de Yang-Baxter}
\end{figure}
\end{center}

Pour d\'ecrire des solutions de l'\'equation de Yang-Baxter quantique, on peut sp\'ecialiser sur des {\bf repr\'esentations} de dimension finie de $\U_q(\hat{\Glie})$.
Une repr\'esentation (lin\'eaire) de $\U_q(\hat{\Glie})$ est un espace vectoriel $V$ (ici complexe) muni d'un morphisme d'alg\`ebre
$$\rho_V : \U_q(\hat{\Glie}) \rightarrow \text{End}(V).$$
Autrement dit, l'alg\`ebre $\U_q(\hat{\Glie})$ agit sur l'espace $V$ par op\'erateurs lin\'eaires. 

L'\'etude des repr\'esentations est un vaste domaine, central en math\'ematiques, appel\'e {\bf th\'eorie des repr\'esentations}. En arithm\'etique par exemple, les repr\'esentations 
de groupes de Galois jouent un r\^ole crucial. Elles sont \'egalement essentielles dans la formulation m\^eme des principes de la physique quantique car ils font intervenir des repr\'esentations
de l'alg\`ebre des observables.

On d\'efinit naturellement la {\bf somme directe de repr\'esentations} $(V,\rho_V)$ et $(V',\rho_{V'})$ avec l'application $\rho_{V\oplus V'} = \rho_V + \rho_{V'}$ \`a valeurs
dans $\text{End}(V \oplus V')$.

Les {\bf repr\'esentations simples}, c'est-\`a-dire qui n'ont pas de sous-repr\'esentation (sous-espace stable pour l'action de l'alg\`ebre) propre, sont particuli\`erement importantes, comme nous allons le voir dans notre \'etude.
Elles constituent les \og briques \'el\'ementaires\fg{} de la th\'eorie des repr\'esentations. Par exemple, toute repr\'esentation de dimension finie de $\U_q(\Glie)$ est {\bf semi-simple}, c'est-\`a-dire isomorphe \`a une somme
directe de repr\'esentations simples\footnote{Ce r\'esultat d\'emontr\'e par M. Rosso et G. Lusztig est un analogue quantique du th\'eor\`eme classique de Weyl qui assure que toute repr\'esentation de dimension finie de $\U(\Glie)$ est semi-simple.}. Ce n'est pas le cas\footnote{Cependant, toute repr\'esentation $V$ de dimension finie de $\U_q(\hat{\Glie})$ admet une filtration de Jordan-H\"older par des sous-repr\'esentations $V_0 = V\supset V_1\supset V_2 \cdots \supset V_N = \{0\}$  avec les  $V_i/V_{i+1}$ simples.} pour l'alg\`ebre affine quantique $\U_q(\hat{\Glie})$.

Comme $\U_q(\hat{\Glie})$ est munie d'un coproduit (\ref{coproduit}), pour deux repr\'esentations $(V,\rho_V)$ et $(V',\rho_{V'})$, le produit tensoriel $V\otimes V'$
est aussi une repr\'esentation en utilisant
$$\rho_{V\otimes V'} = (\rho_V\otimes \rho_{V'})\circ \Delta : \U_q(\hat{\Glie})\rightarrow \text{End}(V)\otimes \text{End}(V') = \text{End}(V\otimes V').$$
Cette action sur un {\bf produit tensoriel de repr\'esentation} sera utile par la suite. Mais ind\'ependamment on peut faire aussi agir directement la $R$-matrice
universelle sur un carr\'e tensoriel $V\otimes V$ pour $V$ une repr\'esentation de dimension finie de $\U_q(\hat{\Glie})$ : on peut en effet consid\'erer l'image de la $R$-matrice universelle dans $\text{End}(V^{\otimes 2})(z)$
$$\mathcal{R}_{V,V}(z) = (\rho_V\otimes \rho_V)(\mathcal{R}(z))\in \text{End}(V)^{\otimes 2}[[z]] = \text{End}(V^{\otimes 2})[[z]].$$
On obtient aussi une solution de l'\'equation de Yang-Baxter quantique, dite {\bf $R$-matrice}, mais dans l'alg\`ebre de dimension finie $\text{End}(V^{\otimes 2})[[z]]$.

Par exemple, dans le cas $\Glie = sl_2$, l'alg\`ebre affine quantique $\U_q(\hat{sl_2})$ poss\`ede une repr\'esentation de dimension $2$ dite {\bf repr\'esentation fondamentale}
et not\'ee $V_1$. Par le proc\'ed\'e d\'ecrit ci-dessus, elle produit la $R$-matrice suivante
\footnote{La solution explicite de l'\'equation de Yang-Baxter donn\'ee ici est la $R$-matrice \og normalis\'ee\fg, obtenue en multipliant $\mathcal{R}_{V_1, V_1}(z)$ par une certaine fonction scalaire de $z$. On peut constater que ses coefficients sont des fonctions m\'eromorphes en $z$. C'est un ph\'enom\`ene g\'en\'eral, voir \cite{efk}.} dans $\text{End}(V_1^{\otimes 2})[[z]]$ avec $V_1^{\otimes 2}$ qui est de dimension $4$ :
$$\begin{pmatrix}1&0&0&0\\ 0&\frac{q^{-1}(z-1)}{z-q^{-2}}&\frac{z(1 - q^{-2})}{z-q^{-2}}&0\\0&\frac{1-q^{-2}}{z - q^{-2}}&\frac{q^{-1}(z - 1)}{z - q^{-2}}&0\\0&0&0&1\end{pmatrix}.$$
C'est la $R$-matrice associ\'ee au mod\`ele $XXZ$. Mais la th\'eorie des groupes quantiques en produit beaucoup d'autres, selon qu'on change l'alg\`ebre de Lie $\Glie$ ou la repr\'esentation $V$. 
Elles correspondent \`a autant de syst\`emes quantiques. 

La {\bf matrice de transfert} $\mathcal{T}_V(z)$ est alors d\'efinie en prenant la trace partielle sur la repr\'esentation, c'est-\`a-dire 
\begin{equation}\label{transfer}
\mathcal{T}_V(z) = ((\on{Tr}_V \circ \rho_V) \otimes \on{id})({\mathcal{R}(z)})\in \U_q(\hat{\Glie})[[z]].
\end{equation}
La repr\'esentation $V$ qui sert \`a construire la matrice de transfert $\mathcal{T}_V(z)$ est appel\'ee {\bf espace auxiliaire}.
Comme cons\'equence de l'\'equation de Baxter, les matrices de transfert commutent, c'est-\`a-dire que pour une autre repr\'esentation $V'$ on a
$$\mathcal{T}_V(z)\mathcal{T}_V(z') = \mathcal{T}_V(z')\mathcal{T}_V(z)\text{ dans }\U_q(\hat{\Glie})[[z,z']].$$
Ainsi, les coefficients $\mathcal{T}_V[N]$ des matrices de transfert, d\'efinis par
$$\mathcal{T}_V(z) = \sum_{N\geq 0}z^N \mathcal{T}_V[N],$$
engendrent une sous-alg\`ebre commutative de $\U_q(\hat{\Glie})$. 

Donnons-nous une autre repr\'esentation de dimension finie $W$ de $\U_q(\hat{\Glie})$, dite {\bf espace des \'etats}. Les coefficients $\mathcal{T}_V[N]$ des matrices de transfert agissent donc sur $W$ en une grande famille commutative d'op\'erateurs. Ainsi, il fait sens de parler des valeurs propres des matrices de transfert $\mathcal{T}_V(z)$ sur $W$.

Dans le cas particulier du mod\`ele $XXZ$, on rappelle que ${\mathfrak g} = sl_2$ et $V = V_1$ est une repr\'esentation fondamentale de dimension $2$.
L'espace des \'etats $W$ est un produit tensoriel de repr\'esentations fondamentales de dimension $2$ et l'image de l'op\'erateur $\mathcal{T}_{V_1}(z)$
dans $\text{End}(W)[[z]]$ est bien la matrice de transfert de Baxter. Les r\'esultats de Baxter donnent donc la structure du spectre de $\mathcal{T}_{V_1}(z)$
sur $W$ dans ce cas. 

\medskip

Que dire en g\'en\'eral ?

\section{La conjecture du spectre quantique}

En 1998 \cite{Fre}, E. Frenkel et N. Reshetikhin ont propos\'e une nouvelle approche dans le but de g\'en\'eraliser les formules de Baxter.

\`A cette fin, ils ont introduit le {\bf $q$-caract\`ere} $\chi_q(V)$ d'une repr\'esentation $V$ de dimension finie de $\U_q(\hat{\Glie})$. Il s'agit d'un polyn\^ome
de Laurent \`a coefficients entiers en des ind\'etermin\'ees $Y_{i,a}$ ($1\leq i\leq n$, $a\in\CC^*$) 
$$\chi_q(V) \in \ZZ[Y_{i,a}^{\pm 1}]_{1\leq i\leq n, a\in\CC^*}.$$
L'entier $n$ est ici le {\bf rang} de l'alg\`ebre de Lie $\Glie$, qui par exemple vaut bien $n$ pour $\Glie = sl_{n+1}$. La d\'efinition du $q$-caract\`ere de $V$ repose
sur une d\'ecomposition de $V$ en sous-espaces de Jordan\footnote{Pour une famille commutative d'op\'erateurs sur $W$, obtenus \`a partir de la r\'ealisation de Drinfeld de $\U_q(\hat{\mathfrak{g}})$ et distincts en g\'en\'eral des coefficients des matrices de transfert.} $V_m$ param\'etr\'es par des mon\^ome $m$ en les variables $Y_{i,a}^{\pm 1}$ :
$$V = \bigoplus_m V_m.$$
Le $q$-caract\`ere encode les dimensions de cette d\'ecomposition. Il est d\'efini par
$$\chi_q(V) = \sum_m \text{dim}(V_m) m.$$
Ainsi, les coefficients de $\chi_q(V)$ sont en fait positifs et leur somme est la dimension $V$. 

Par exemple, pour $\Glie = sl_2$ et $V = V_1$ la repr\'esentation fondamentale de dimension $2$,  
\begin{equation}\label{qcar}
\chi_q(V) = Y_{1,q^{-1}} + Y_{1,q}^{-1}.
\end{equation}
On a donc dans ce cas deux sous-espaces de Jordan de dimension $1$ associ\'es aux mon\^omes respectifs $Y_{1,q^{-1}}$ et $Y_{1,q}^{-1}$ :
$$V = V_{Y_{1,q^{-1}}} \oplus V_{Y_{1,q}^{-1}}.$$

La {\bf conjecture du spectre quantique} de Frenkel et Reshetikhin \cite{Fre} pr\'edit\footnote{Dans des cas particuliers, une conjecture analogue avait \'et\'e formul\'ee par N. Reshetikhin \cite{R3}; V. Bazhanov et N. Reshetikhin
\cite{BR}; et A. Kuniba et J. Suzuki \cite{KS}.} que pour une repr\'esentation de dimension finie donn\'ee $V$,
les valeurs propres $\lambda_j$ de $\mathcal{T}_V(z)$ sur une repr\'esentation simple\footnote{
Plus g\'en\'eralement, $W$ peut \^etre un produit tensoriel de repr\'esentations simples.} $W$ sont obtenues de la mani\`ere suivante :
dans le $q$-caract\`ere $\chi_q(V)$ de $V$, on remplace chaque variable formelle $Y_{i,a}$ par\footnote{Pour simplifier l'exposition, on supposera dans la suite de $\Glie$
est simplement lac\'ee (c'est le cas notamment des alg\`ebres de Lie $sl_{n+1}$).}
$$F_{i}(az) q^{\text{deg}(Q_{i,j})} \frac{Q_{i,j}(zaq^{-1})}{Q_{i,j}(zaq)},$$
o\`u $F_{i}(z)$ est une fonction universelle, au sens o\`u elle ne d\'epend pas de la valeur propre $\lambda_j$, et $Q_{i,j}(z)$ d\'epend de la valeur propre $\lambda_j$ mais est un polyn\^ome. C'est l'analogue du polyn\^ome de Baxter.

Notons que c'est bien le $q$-carat\`ere {\it de l'espace auxiliaire} $V$ qui est utilis\'e pour \'ecrire la formule du spectre de la matrice de transfert {\it sur l'espace des \'etats} $W$.

Dans le cas particulier du mod\`ele $XXZ$, on obtient  \`a partir de (\ref{qcar}) la formule
$$\lambda_j = F_{1}(zq^{-1}) q^{\text{deg}(Q_{1,j})} \frac{Q_{1,j}(zq^{-2})}{Q_{1,j}(z)} + (F_1(zq))^{-1} q^{-\text{deg}(Q_{1,j})} \frac{Q_{1,j}(zq^2)}{Q_{1,j}(z)}.$$
Ainsi, la conjecture est bien compatible avec la formule de Baxter (\ref{relB}) en identifiant 
$$A(z) = (D(zq^2))^{-1} = (F_1(zq))^{-1}q^{-\text{deg}(Q_{1,j})}.$$ 
On peut d\'etailler par exemple le cas o\`u l'espace des \'etats $W\simeq V_1$ est de dimension $2$. On a alors $2$ valeurs propres $\lambda_0$ et $\lambda_1$. La fonction universelle est
$$F_1(z) = q^{1/2}\text{exp}\left(  \sum_{r > 0}  \frac{z^r (q^{-r} - q^r)}{r(q^r + q^{-r})}\right),$$
et les polyn\^omes de Baxter sont
$$Q_{1,0}(z) = 1\text{ et }Q_{1,1}(z) = 1 - z(1 + q + q^2).$$
On obtient donc le spectre
$$\lambda_0 = F_1(zq^{-1})\left(1 + q^{-3}\frac{1 - z^{-1}}{1 - z^{-1}q^{-2}}\right),$$
\begin{equation}\label{bethe}\lambda_1 = F_1(zq^{-1}) \left(q\frac{1-z(1+q^{-1} + q^{-2})}{1 - z(1 + q + q^2)} + q^{-4}\frac{(1 - z^{-1})(1 - z(q^2 + q^3+ q^4))}{(1 - z^{-1}q^{-2})(1 - z(1 + q + q^2))}\right).\end{equation}

En g\'en\'eral la formule peut avoir plus de deux termes.
Par exemple, dans le cas d'une certaine repr\'esentation fondamentale $V$ de dimension $3$ de $\U_q(\hat{sl_3})$, le $q$-caract\`ere est
\begin{equation}\label{sl3}\chi_q(V) = Y_{1,q^{-1}} + Y_{1,q}^{-1}Y_{2,1} + Y_{2,q^2}^{-1},\end{equation}
et la formule pour le spectre est
$$F_1(zq^{-1}) q^{\text{deg}(Q_{1,j})}\frac{Q_{1,j}(zq^{-2})}{Q_{1,j}(z)} +\frac{F_2(z) q^{\text{deg}(Q_{2,j})}}{F_1(zq)q^{\text{deg}(Q_{1,j})}} \frac{Q_{1,j}(zq^2)Q_{2,j}(zq^{-1})}{Q_{1,j}(z)Q_{2,j}(zq)} + \frac{q^{-\text{deg}(Q_{2,j})}}{(F_2(zq^2))^{-1}} \frac{Q_{2,j}(zq^3)}{Q_{2,j}(zq)}.$$

Notons qu'en g\'en\'eral les repr\'esentations simples $V$ de dimension finie peuvent avoir une dimension \og tr\`es grande\fg. Par exemple, H. Nakajima a obtenu (\`a l'aide d'un super-calculateur et en s'appuyant sur \cite{Nak}) que dans le cas de l'alg\`ebre de Lie exceptionelle de type $E_8$, une des repr\'esentations fondamentales a un $q$-caract\`ere avec $6899079264$ mon\^omes qui n\'ecessite un fichier de taille m\'emoire $180$ Go pour \^etre \'ecrit. Il y a donc autant de termes dans la formule de Baxter correspondante. Et les repr\'esentations fondamentales sont les repr\'esentations simples de dimensions les plus basses. 

Il est donc hors de question d'aborder cette conjecture par un calcul explicite en g\'en\'eral. D'ailleurs, m\^eme si les repr\'esentations simples de dimension finie de $\U_q(\hat{\Glie})$ ont \'et\'e intensivement \'etudi\'ees ces vingt-cinq derni\`eres ann\'ees, on ne conna\^it pas en g\'en\'eral de formule pour leur $q$-caract\`ere, ni m\^eme en fait pour leur dimension.

Ainsi, il faut de nouvelles structures pour aborder la conjecture du spectre quantique.

Notre d\'emonstration avec E. Frenkel \cite{FH} de la conjecture du spectre quantique repose ainsi sur 
de nouveaux ingr\'edients dont nous donnons un bref aper\c cu dans les sections suivantes.

\section{Repr\'esentations pr\'efondamentales}

L'id\'ee g\'en\'erale de la preuve est d'interpr\'eter les $Q_i$ eux-m\^emes comme des valeurs
propres de nouvelles matrices de transfert, construites non pas \`a partir de repr\'esentations de dimension
finie $V$, mais de repr\'esentations de dimension infinie dite {\bf repr\'esentations pr\'efondamentales}
$L_{i,a}^+$ o\`u $1\leq i\leq n$ et $a\in\CC^*$.

Nous avions construit pr\'ealablement ces repr\'esentations pr\'efondamentales avec M. Jimbo \cite{HJ} dans un contexte un peu diff\'erent. Ce ne sont pas des repr\'esentations de l'alg\`ebre enti\`ere $\U_q({\hat{\mathfrak g}})$, mais d'une certaine sous-alg\`ebre, la {\bf sous-alg\`ebre de Borel} 
$$\U_q(\hat{\mathfrak{b}})\subset \U_q({\hat{\mathfrak g}}).$$ 
Cela ne pose cependant pas de probl\`eme pour construire la matrice de transfert $\mathcal{T}_{i,a}(z)$ associ\'ee \`a la repr\'esentation pr\'efondamentale $L_{i,a}$ par la formule (\ref{transfer}), car il est justement connu que la partie \og gauche\fg{} de la $R$-matrice universelle (celle \`a qui on applique $\rho_{L_{i,a}}$) est dans la sous-alg\`ebre de Borel\footnote{On ne peut cependant pas appliquer la trace \`a un espace de dimension infinie. On utilise une graduation naturelle de $L_i$ par des espaces de dimension finie (les espaces de poids). Ainsi dans la suite, les traces, matrices de transfert, etc. sont \og tordues\fg{} par cette graduation.} :
$$\mathcal{T}_{i,a}(z) =  ((\on{Tr}_{L_{i,a}} \circ \rho_{L_{i,a}}) \otimes \on{id})({\mathcal{R}(z)})\in \U_q(\hat{\Glie})[[z]].$$
Il n'est alors pas difficile de montrer qu'en utilisant un certain automorphisme de $\U_q(\hat{\mathfrak{b}})$ on a 
$$\mathcal{T}_{i,a}(z) = \mathcal{T}_i(az)\text{ o\`u }\mathcal{T}_i(z) = \mathcal{T}_{i,1}(z).$$

Pour le cas du mod\`ele $XXZ$, c'est-\`a-dire pour $\Glie = sl_2$, V. Bazhanov, S. Lukyanov, et A. Zamolodchikov avaient d\'ej\`a construit \og \`a la main\fg{} une repr\'esentation pr\'efondamentale (appel\'ee repr\'esentation de $q$-oscillation) et la matrice de transfert associ\'ee (appel\'ee $Q$-op\'erateur de Baxter) dans l'article important \cite{BLZ}. 

Pour obtenir l'existence des repr\'esentations pr\'efondamentales en g\'en\'eral \cite{HJ}, on ne peut encore une fois pas faire de calculs explicites : le point crucial est de consid\'erer des syst\`emes inductifs\footnote{Les inclusions $L_k\subset L_{k+1}$, construites \`a l'aide de produits tensoriels de sous-espaces \cite{h3}, ne sont pas compatibles avec l'action de $\U_q(\hat{\Glie})$ enti\`ere mais avec celle d'une sous-alg\`ebre $\U_q^+(\hat{\mathfrak{b}})$ de $\U_q(\hat{\mathfrak{b}})$.} de repr\'esentations simples $L_k$ (les repr\'esentations de Kirillov-Reshetikhin) de dimension finie strictement croissante avec $k\geq 0$ et de d\'eterminer en quel sens l'action de la sous-alg\`ebre de Borel $\U_q(\hat{\mathfrak{b}})$ \og converge\fg{} sur la limite inductive $L_\infty$, qui elle est de dimension infinie :
$$L_0\subset L_1\subset L_2\subset \cdots \subset L_k\subset L_{k+1}\subset \cdots  \subset L_{\infty}.$$
Il s'agit ainsi d'une construction asymptotique des repr\'esentations pr\'efondamentales.

En utilisant certaines filtrations de la repr\'esentation pr\'efondamentale $L_{i,a}$, nous \'etablissons qu'effectivement, \`a un facteur scalaire universel $f_i(z)$ pr\`es, la matrice de transfert associ\'ee $\mathcal{T}_i(z)$ agit sur l'espace des \'etats $W$ par un op\'erateur polyn\^omial :
$$\rho_W(\mathcal{T}_i(z)) \in f_i(z) \times (\text{End}(W))[z].$$
Il n'est pas difficile d'\'ecrire une formule explicite pour la fonction universelle scalaire $f_i(z)\in\CC[[z]]$ (elle ne d\'epend que de $V$ et de $W$). Il est beaucoup plus d\'elicat d'obtenir des informations sur la partie lin\'eaire polyn\^omiale
$$(f_i(z))^{-1}\rho_W(\mathcal{T}_i(z)) \in (\text{End}(W))[z].$$

De m\^eme que les matrices de transfert usuelles commutent, on a 
$$\mathcal{T}_i(z)\mathcal{T}_i(z') = \mathcal{T}_i(z')\mathcal{T}_i(z),$$
et donc on obtient une famille commutative $\mathcal{T}_i[m]$ si on \'ecrit
$$\mathcal{T}_i(z) = \sum_{m\geq 0}\mathcal{T}_i[m]z^m.$$
En utilisant la trigonalisation simultan\'ee, cette commutativit\'e implique que les valeurs propres sur $W$ de $(F_i(z))^{-1}\mathcal{T}_i(z)$ elles-m\^emes sont \'egalement des polyn\^omes.

\section{Anneau de Grothendieck et relations de Baxter}

Il faut enfin d\'emontrer que les valeurs propres de la matrice de transfert $\mathcal{T}_V(z)$
s'expriment, comme pr\'evu dans la conjecture, en terme des valeurs propres des $\mathcal{T}_i(z)$
selon le $q$-caract\`ere de $V$. Autrement dit, en rempla\c cant dans $\chi_q(V)$ chaque
variable $Y_{i,a}$ par le quotient\footnote{Ce quotient doit en fait \^etre multipli\'e par une matrice de transfert d'une repr\'esentation de dimension $1$ que nous omettons dans la suite pour simplifier l'exposition.} 
$$\mathcal{T}_i(azq^{-1})/\mathcal{T}_i(azq),$$ 
obtient-on la matrice de transfert $\mathcal{T}_V(z)$ ?

Dans le cas ${\mathfrak g} = \wh{sl}_2$ et $V$ de dimension du mod\`ele $XXZ$, un calcul \cite{BLZ} donne le r\'esultat. On a bien : 
$$\mathcal{T}_V(z) = \frac{\mathcal{T}_1(zq^{-1})}{\mathcal{T}_1(zq)} + \frac{\mathcal{T}_1(zq^3)}{\mathcal{T}_1(zq)}.$$

En g\'en\'eral, nous proposons d'utiliser la {\bf cat\'egorie $\mathcal{O}$} que nous avons d\'efinie avec M. Jimbo \cite{HJ}. Il s'agit d'une {\bf cat\'egorie mono\" idale} (stable par produits tensoriels) de repr\'esentations de l'alg\`ebre de Borel $\U_q(\hat{\mathfrak{b}})$, contenant les repr\'esentations de dimension finie ainsi que les repr\'esentations pr\'efondamentales. 
Nous {\bf cat\'egorifions} les relations de Baxter g\'en\'eralis\'ees, c'est-\`a-dire que nous les exprimons en termes de la cat\'egorie $\mathcal{O}$. Pour ce faire, on peut d\'efinir l'{\bf anneau de Grothendieck} $K(\mathcal{O})$ de cette cat\'egorie. En tant que groupe, il s'agit du groupe libre engendr\'e par les classes d'isomorphismes de repr\'esentations simples :
$$K(\mathcal{O}) = \bigoplus_{[V]\text{ Classe d'un simple dans }\mathcal{O}.} \ZZ [V].$$
Alors tout objet (non n\'ecessairement simple) de $\mathcal{O}$ a une image dans $K(\mathcal{O})$ en imposant la relation
$$[V''] = [V] + [V']$$
si on a une suite exacte dans la cat\'egorie
$$0\rightarrow V \rightarrow V''\rightarrow V'\rightarrow 0.$$
On peut alors munir $K(\mathcal{O})$ d'une structure d'anneau par la relation
$$[V\otimes V'] = [V][V']$$
pour des objets $V$, $V'$ de la cat\'egorie $\mathcal{O}$.

Un des th\'eor\`emes principaux de \cite{FH} est qu'en rempla\c cant dans $\chi_q(V)$ chaque
variable $Y_{i,a}$ par le quotient 
$$\frac{[L_{i,aq^{-1}}]}{[L_{i,aq}]},$$ 
en repla\c cant $\chi_q(V)$ par $[V]$ puis en \og chassant\fg{} les d\'enominateurs, on obtient une relation dans l'anneau de Grothendieck $K(\mathcal{O})$.

Par exemple, dans notre cas favori du mod\`ele $XXZ$, on obtient
$$[V] = \frac{[L_{1,q^{-1}}]}{[L_{1,q}]} + \frac{[L_{1,q^3}]}{[L_{1,q}]}$$
qui donne la relation de Baxter cat\'egorifi\'ee dans l'anneau de Grothendieck
$$[V][L_{1,q}] = [V\otimes L_{1,q}] = [L_{1,q^{-1}}] + [L_{1,q^3}].$$
En g\'en\'eral on obtient des relations avec plus de termes, comme dans l'exemple pour $\Glie = sl_3$ ci-dessus pour lequel la formule (\ref{sl3}) donne
$$[V\otimes L_{1,1}\otimes L_{2,q}] = [L_{1,q^{-2}}\otimes L_{2,q}] + [L_{1,q^2}\otimes L_{2,q^{-1}}] + [L_{2,q^3}\otimes L_{1,1}].$$
Maintenant, \og prendre la matrice de transfert\fg{} est additif et multiplicatif, c'est \`a dire qu'on a un morphisme d'anneau\footnote{On peut montrer que ce morphisme d'anneau est injectif et donc que l'anneau de Grothendieck $K(\mathcal{O})$ est commutatif (bien que la cat\'egorie ne soit pas tress\'ee, c'est-\`a-dire que $V\otimes V'$ et $V'\otimes V$ ne sont pas isomorphes en g\'en\'eral). En fait, l'application de $q$-caract\`eres $[V]\mapsto \chi_q(V)$ elle-m\^eme peut \^etre prolong\'ee en un morphisme d'anneau injectif sur $K(\mathcal{O})$.}
$$\mathcal{T} : K(\mathcal{O})\rightarrow \mathcal{U}_q(\hat{\Glie})[[z]]\text{ , }[V]\mapsto \mathcal{T}_V(z).$$
Ainsi, les relations de Baxter g\'en\'eralis\'ees dans l'anneau de Grothendieck $K(\mathcal{O})$ impliquent les relations voulues entre les matrices de transfert. La conjecture du spectre quantique est donc d\'emontr\'ee.

\begin{center}*\end{center}

Pour conclure, les formules pour les valeurs propres des matrices de transfert en terme
des polyn\^omes $Q_{i,j}$ impliquent des \'equations entre les racines de ces polyn\^omes pour garantir que les p\^oles
apparents se simplifient (par exemple dans l'\'equation (\ref{bethe}), $(1 + q + q^2)^{-1}$ n'est en fait pas un p\^ole de $\lambda_1$). Dans le cas du mod\`ele $XXZ$ ce sont les fameuses \'equations de l'Ansatz de Bethe.
Ces consid\'erations ont men\'e N. Reshetikhin \cite{R3} \`a formuler ces \'equations dans le cas g\'en\'eral
  (voir aussi \cite{BR,KS, F:icmp}). La preuve de la conjecture du spectre quantique permet de donner une explication et une approche uniforme \`a ces formules. On a maintenant une autre conjecture importante et ouverte : l'existence d'une bijection entre toutes les valeurs propres et les solutions des \'equations de l'Ansatz de Bethe (conjecture de compl\'etude).

\medskip

{\bf Remerciements} : je souhaite adresser mes remerciements \`a E. Ghys pour m'avoir encourag\'e \`a \'ecrire cet article, \`a E. Frenkel et M. Jimbo pour notre collaboration et enfin \`a J. Dumont, C. Hernandez, P. Zinn-Justin et l'\'equipe d'Images de Math\'ematiques pour leurs remarques sur une version pr\'eliminaire de ce texte.

\backmatter


\begin{thebibliography}{ASM}

\bibitem[B1]{Baxter} {R.J. Baxter},
\newblock Partition Function of the Eight-Vertex Lattice Model,
\newblock Ann. Phys. {\bf 70} (1971), 193--228.

\bibitem[B2]{BaxterBook} {R.J. Baxter},
\newblock Exactly solved models in statistical mechanics, 
\newblock Academic Press Inc., London (1982).

\bibitem[BLZ]{BLZ} {V.V. Bazhanov, S.L. Lukyanov et
  A.B. Zamolodchikov},
\newblock Integrable structure of conformal field
theory. II. Q-operator and DDV equation,
\newblock Comm. Math. Phys. {\bf 190} (1997), 247--278.

\bibitem[BR]{BR} {V.V. Bazhanov et N.Yu. Reshetikhin},
\newblock Restricted solid on solid models connected with simply laced
Lie algebra,
\newblock J. Phys. {\bf A 23} (1990), 1477--1492.

\bibitem[CP]{CP} {V. Chari et A. Pressley},
\newblock Guide to Quantum Groups,
\newblock Cambridge University Press, Cambridge, 1994.

\bibitem[D1]{Dri} V. Drinfel'd, 
\newblock Quantum groups,
\newblock Proceedings of the International Congress of Mathematicians, Vol. 1, 2 (Berkeley, Calif., 1986), 798-820, Amer. Math. Soc., Providence, RI, (1987)

\bibitem[D2]{Dri2} {V. Drinfel'd}, 
\newblock A new realization of Yangians and of 
quantum affine algebras, 
\newblock Soviet Math. Dokl. {\bf 36} (1988), 212--216.  

\bibitem[EFK]{efk} {P. Etingof, I. Frenkel and A. Kirillov},
\newblock Lectures on representation theory and Knizhnik-Zamolodchikov
equations,
\newblock Mathematical Surveys and Monographs {\bf 58}. American Mathematical Society, Providence, RI, 1998.

\bibitem[F]{F:icmp} {E. Frenkel},
\newblock Affine algebras, Langlands duality and Bethe ansatz,
\newblock in Proceedings of the International Congress of
Mathematical Physics, Paris, 1994, ed. D. Iagolnitzer, pp. 606--642,
International Press, 1995.

\bibitem[FH]{FH} {E. Frenkel et D. Hernandez},
\newblock Baxter's Relations and Spectra of Quantum Integrable Models,
\newblock Pr\'epublication arXiv:1308.3444, \`a para\^itre dans Duke. Math. J.

\bibitem[FR]{Fre} {E. Frenkel et N. Reshetikhin}, 
\newblock The $q$-characters of representations of quantum 
affine algebras and deformations of $W$-Algebras, 
\newblock in Recent Developments in Quantum Affine Algebras and 
related topics, 
\newblock 
Contemp. Math. {\bf 248} (1999), 163--205.

\bibitem[H]{h3} {D. Hernandez}, 
\newblock Simple tensor products, 
\newblock Invent. Math. {\bf 181} (2010), 649--675.

\bibitem[HJ]{HJ} {D. Hernandez et M. Jimbo},
\newblock Asymptotic representations and Drinfeld rational fractions,
\newblock Compos. Math. {\bf 148} (2012), no. 5, 1593--1623.

\bibitem[J]{J} M. Jimbo, 
\newblock A q-difference analogue of $\mathcal{U}(\mathfrak{g})$ and the Yang-Baxter equation, 
\newblock Lett. Math. Phys. {\bf 10} (1985), no. 1, 63-69.

\bibitem[JM]{JM} M. Jimbo et T. Miwa, 
\newblock Algebraic analysis of solvable lattice models, 
\newblock CBMS Regional Conference Series in Mathematics {\bf 85}, American Mathematical Society (1995).

\bibitem[K]{ka}{V. Kac},
\newblock Infinite Dimensional Lie Algebras, 3rd ed.,
\newblock Cambridge University Press, Cambridge, 1990.

\bibitem[KS]{KS} {A. Kuniba et J. Suzuki},
\newblock Analytic Bethe Ansatz for Fundamental Representations of
Yangians,
\newblock Comm. Math. Phys. {\bf 173} (1995), 225---264.

\bibitem[N]{Nak} {H. Nakajima},
\newblock Quiver varieties and t-analogs of q-characters of quantum affine algebras,
\newblock Ann. of Math. (2) {\bf 160} (2004), no. 3, 1057--1097.

\bibitem[R]{R3} {N. Reshetikhin,}
\newblock The spectrum of the transfer matrices connected with
Kac--Moody algebras,
\newblock Lett. Math. Phys. {\bf 14} (1987), 235--246.

\end{thebibliography}
\end{document}